\documentclass[12pt]{article}

\newfam\msbfam
\font\tenmsb=msbm10  scaled \magstep1 \textfont\msbfam=\tenmsb
\font\sevenmsb=msbm7  scaled \magstep1 \scriptfont\msbfam=\sevenmsb
\font\fivemsb=msbm5  scaled \magstep1 \scriptscriptfont\msbfam=\fivemsb
\def\Bbb{\fam\msbfam \tenmsb}

\def\RR{{\Bbb R}}
\def\CC{{\Bbb C}}

\newtheorem{theorem}{Theorem}

\newtheorem{lemma}[theorem]{Lemma}

\newtheorem{definition}{Definition}
\newtheorem{example}[definition]{EXAMPLE}

\def\ss{\subseteq}

\def\ra{\rightarrow}
\def\O{\Omega}

 \def\HollowBox #1#2{{\dimen0=#1 \advance\dimen0 by -#2       
       \dimen1=#1 \advance\dimen1 by #2                       
        \vrule height #1 depth #2 width #2                    
        \vrule height 0pt depth #2 width #1                   
        \llap{\vrule height #1 depth -\dimen0 width \dimen1}%
       \hskip -#2                                             
       \vrule height #1 depth #2 width #2}}                   
 \def\BoxOpTwo{\mathord{\HollowBox{6pt}{.4pt}}\;}             

\def\endpf{\hfill $\BoxOpTwo$}

\begin{document}

\begin{center}
{\LARGE \bf Topologies on the Space}
\medskip \\
{\LARGE \bf of Holomorphic Functions (REVISED)}
\end{center}
\vspace*{.2in}

\begin{center}
Steven G. Krantz\footnote{The author thanks the American Institute of Mathematics
for its hospitality and support during a portion of this work.}
\end{center}

\section{Introduction}

{\it We remark that the result presented here can be proved more simply using
the closed graph theorem.  However we believe that our results can be used
to prove a more interesting result.  Details to follow in a later paper.}

One of the remarkable features of the space of holomorphic functions
(in either one or several variables) is that the standard Frechet space
topologies---say, for example, the $L^2$ or Bergman norm---control a stronger
(and simpler) and more useful topology, namely uniform convergence on compact sets.
This simple fact lies at the heart of many key results in basic
complex function theory---for example the completeness of many
important function spaces.

It is natural to wonder whether this property is universal.  Is it
the case that {\it any} Frechet space topology on the space of holomorphic
functions implies uniform convergence on compact sets (equivalently, convergence in the
compact-open topology)?  The surprising answer to this question---suitably
formulated---is ``yes'', and that is the main result of the present paper.

It is a pleasure to thank Peter Pflug for posing the question that led to
this paper, and for early discussions of the matter.

\section{Basic Concepts}

A {\it domain} is a connected open set.   We typically denote a domain
with the symbol $\Omega$.  Let ${\cal O}(\Omega)$ denote the complex linear
space of holomorphic functions on $\Omega$.

A {\it Frechet space topology}\ ${\cal T}$ on ${\cal O}(\Omega)$ is a collection
$\{{\cal F}_\alpha\}_{\alpha \in A}$ of semi-norms on ${\cal O}(\O)$ so that
the resulting topology is complete.   Many of the most useful topologies on
${\cal O}(\Omega)$ are Frechet space topologies; some, however, are Banach
space topologies.

Our basic result is this:

\begin{theorem} \sl
Let $\O \ss \CC^n$ be a pseudoconvex domain.  Let ${\cal T}$ be a Frechet space topology on ${\cal O}(\O)$.
Assume that each semi-norm ${\cal F}_\alpha$ in ${\cal T}$ is majorized by a semi-norm ${\cal G}_\alpha$
on $C^\infty(\Omega)$.  If a sequence of functions $\{f_j\} \ss {\cal O}(\O)$ converges in 
the topology ${\cal T}$, then it converges uniformly on compact sets.
\end{theorem}

The main tool in the proof of the theorem is a fairly general version of the Hahn-Banach
theorem, for which see [YOS, p.\ 105, ff.]:

\begin{theorem} \sl
Let $V$ be a vector space over the field $\CC$.  A function $N: V \ra \RR$ is called
{\it sublinear} if
$$
N(ax + by) \leq |a| N(x) + |b| N(y)
$$
for all $x, y \in V$ and all complex scalars $a, b$.  Let $U$ be a complex linear subspace of $V$ 
(not necessarily closed) and let $\varphi: U \ra \CC$ be a linear functional (not necessarily
bounded) on $U$.  If $\varphi$ is dominated by $N$ in the sense that $|\varphi(x)| \leq N(x)$ for all $x \in U$, then there is
a linear extension $\widehat{\varphi}: V \ra \CC$ of $\varphi$ to all of $V$ (meaning that
$\widehat{\varphi} \bigr |_U = \varphi$) which is also dominated by $N$.
\end{theorem}

Since the next ideas do not seem to be well-documented in the literature, we close this
section by reviewing a few elementary concepts from Frechet space theory.

If ${\cal F}_\alpha$ is a semi-norm on our Frechet space $X$, and if $\lambda$ is a linear
functional on $X$, then we may define the {\it dual norm} of $\lambda$ by
$$
{\cal F}^*_\alpha(\lambda) = \sup_{f \in X \atop
                         |{\cal F}_\alpha(f)| \leq 1} | \langle f, \lambda \rangle | \, .
$$
This definition is of course analogous to the norm of a functional on a Banach space.

With this language and notation, we have the following result:

\begin{lemma} \sl
If $f \in X$ and ${\cal F}_\alpha$ is a semi-norm on $X$, then
$$
{\cal F}_\alpha(f) = \sup_{\lambda \in X^* \atop
                              {\cal F}^*_\alpha(\lambda) \leq 1} |\langle f, \lambda \rangle | \, .
$$
\end{lemma}
{\bf Proof:}  The proof follows classical lines.  See [ZEI, p. 6].

\section{Proof of the Main Result}

Let $\lambda$ denote a typical element of the dual space of ${\cal O}(\O)$ equipped with the
Frechet space topology ${\cal T}$.  According to Lemma 3, we must examine the expressions $\tau_\lambda: f \mapsto \langle f, \lambda \rangle$.  Such a function
is of course a linear functional on ${\cal O}(\O)$.  And it is bounded with respect to the topology
of ${\cal F}_\alpha$, any $\alpha$.  Thus it is dominated (in the sense
defined above) by some semi-norm ${\cal G}_\alpha$ on $C^\infty(\Omega)$.  Thus the
Hahn-Banach theorem applies and there is an extension $\widehat{\lambda}$ of $\lambda$
to $C^\infty(\Omega)$ which is also dominated by ${\cal G}_\alpha$.  It follows that the extension
$\widehat{\lambda}$ is a {\it distribution}.

We then know that $\widehat{\lambda}$ is a finite sum of derivatives of measures.  But the Cauchy
estimates tell us that, when the derivative of a measure is acting on holomorphic functions, 
the action is dominated by the measure itself acting on holomorphic functions.  So we may as well
take it that $\widehat{\lambda}$ is integration against a measure $\mu_\alpha$ on $\Omega$.

Thus for each $\alpha$ we have associated to the functional $\lambda$ a measure $\mu_\alpha$.
Now we must think about the support $K_\alpha$ of each $\mu_\alpha$.  Clearly the hull
of holomorphy of $\cup_\alpha K_\alpha$ must be all of $\Omega$, otherwise the original
Frechet space topology on ${\cal O}(\O)$ will not be complete.   But that says immediately
that the topology majorizes uniform convergence on compact sets.  That is the desired conclusion.
\bigskip \bigskip \\
				  
\section{Examples}

We present two examples that show what the theorem says.

\begin{example} \rm
Consider the space ${\cal O}(D(0,2))$ of holomorphic functions on the disc of radius 2,
and consider the norm
$$
\|f\| = \oint_{|\zeta| = 1} |f(\zeta)| \, d\zeta \, .
$$
The theorem says that convergence in $\| \ \ \|$ implies uniform convergence
on compact sets.  This is a way to re-discover the maximum principle.
\endpf
\end{example}

\begin{example} \rm
Consider the space ${\cal O}(D(0,2))$ as in the last example.  Let $\{a_j\}$ be a sequence
in the disc that has a limit point $p$ in the disc.  Then $\ell^2(\{a_j\})$ is a norm
on ${\cal O}(D(0,2))$.  Thus convergence in this norm implies uniform convergence on
compact sets.  That is a non-obvious fact.
\end{example}

\section{Concluding Remarks}

Not surprisingly, the key fact about holomorphic functions that is used in the proof
of the main result here is the Cauchy estimates.  The Cauchy estimates are logically
equivalent to the analyticity of holomorphic functions, so they are at the heart
of the subject.  It would be interesting to identify other natural spaces of
functions for which results like these hold.

It has been noted by Thomas Krainer and Nigel Kalton that the result presented here
can be proved using the closed graph theorem.  But the proof that we give seems
to have some independent interest.
\vfill
\eject

\noindent {\Large \sc References}

\begin{enumerate}

\item[{\bf [YOS]}]  K. Yosida, {\it Functional Analysis}, $6^{\rm th}$ ed., Springer-Verlag,
Berlin, 1980. 

\item[{\bf [ZEI]}]  E. Zeidler, {\it Applied Functional Analysis:  Main Principles
and their Applications}, Springer, New York, 1995.

\end{enumerate}

\end{document}